\documentclass[11pt]{amsart}
\usepackage{amsxtra,tikz-cd,color}
\usepackage{amssymb}
\usepackage{mathtools}
\theoremstyle{plain}
\newcommand{\cleqn}{\setcounter{equation}{0}}
\newcommand{\clth}{\setcounter{theorem}{0}}
\newcommand {\sectionnew}[1]{\section{#1}\cleqn\clth}
\newtheorem{theorem}{Theorem}[section]
\newtheorem{lemma}[theorem]{Lemma}
\newtheorem{definition-theorem}[theorem]{Definition-Theorem}
\newtheorem{proposition}[theorem]{Proposition}
\newtheorem{corollary}[theorem]{Corollary}
\newtheorem{definition}[theorem]{Definition}
\newtheorem{example}[theorem]{Example}
\newtheorem{remark}[theorem]{Remark}
\newtheorem{conjecture}[theorem]{Conjecture}

\newcommand \bth[1] { \begin{theorem}\label{t#1} }
\newcommand \ble[1] { \begin{lemma}\label{l#1} }

\newcommand \bpr[1] { \begin{proposition}\label{p#1} }
\newcommand \bco[1] { \begin{corollary}\label{c#1} }
\newcommand \bde[1] { \begin{definition}\label{d#1}\rm }
\newcommand \bex[1] { \begin{example}\label{e#1}\rm }
\newcommand \bre[1] { \begin{remark}\label{r#1}\rm }
\newcommand \bcj[1] { \begin{conjecture}\label{j#1}\rm }

\renewcommand {\eth} { \end{theorem} }
\newcommand {\ele} { \end{lemma} }

\newcommand {\epr} { \end{proposition} }
\newcommand {\eco} { \end{corollary} }
\newcommand {\ede} { \end{definition} }
\newcommand {\eex} { \end{example} }
\newcommand {\ere} { \end{remark} }
\newcommand {\ecj} { \end{conjecture} }

\newcommand {\enota} { \end{notation} }
\newcommand \thref[1]{Theorem \ref{t#1}}

\newcommand \prref[1]{Proposition \ref{p#1}}
\newcommand \coref[1]{Corollary \ref{c#1}}

\newcommand \deref[1]{Definition \ref{d#1}}

\newcommand \reref[1]{Remark \ref{r#1}}

\def \Rset {{\mathbb R}}         %mathsets
\def \Cset {{\mathbb C}}

\def \Zset {{\mathbb Z}}

\def \Qset {{\mathbb Q}}

\def \B  {{\mathcal{B}}}
\def \CC {\mathcal{C}}
\def \FF {{\mathcal{F}}}

\def \PP {{\mathcal{P}}}

\def \UU {{\mathcal{U}}}

\def \LL {{\mathcal{L}}}
\def \TT {{\mathcal{T}}}
\def \ZZ {{\mathcal{Z}}}
\def \OO {{\mathcal{O}}}

\def \ZZ {{\mathcal{Z}}}

\def \La {\Lambda}

\def \Om {\Omega}

\def \ep {\varepsilon}

\def \lra {\longrightarrow}

\def \hra {\hookrightarrow}

\def \ol {\overline}

\def \id { {\mathrm{id}} }

\def \sign { {\mathrm{sign}} }

\def \g  {\mathfrak{g}}   % Lie algebra letters

\def \n  {\mathfrak{n}}
\def \mm  {\mathfrak{m}}

\newcommand{\kk}{\Bbbk}

\DeclareMathOperator \Der { {\mathrm{Der}} }

\DeclareMathOperator \res {{\mathrm{res}}}
\DeclareMathOperator \charr { {\mathrm{char}} }

\DeclareMathOperator \maxspec { {\mathrm{MaxSpec}}}

\DeclareMathOperator \tr { {\mathrm{tr}} }

\DeclareMathOperator \End { {\mathrm{End}} }

\DeclareMathOperator \red  { {\mathrm{red}}}
\DeclareMathOperator \reg  { {\mathrm{reg}}}
\DeclareMathOperator \spp { {\mathrm{sp}}}

\newcommand \ex {{\bf{ex}}}

\def \B  {{\widetilde{B}}}
\def \A {{\boldsymbol{\mathsf A}}} 
\def \C {{\boldsymbol{\mathsf C}}} 
\def \U {{\boldsymbol{\mathsf U}}} 
\begin{document}
%%%%%%%%%%%%%%%%%%%%%%%%%%%%%%%%%%%%%%%%%%%%%%%%%%%%%%%%%%%%%%%%%%%%%%%%%%%
%%%%%%%%%%%%%%%%%%%%%%    Title    %%%%%%%%%%%%%%%%%%%%%%%%%%%%%%%%%%%%%%%%
\title
{Poisson trace orders}
\author[K. A. Brown]{K. A. Brown}
\address{School of Mathematics and Statistics \\
University of Glasgow \\
Glasgow G12 8QQ, Scotland}
\email{Ken.Brown@glasgow.ac.uk}
\author[M. T. Yakimov]{M. T. Yakimov}
\thanks{The research of the first named author is supported by Leverhulme Emeritus Fellowship EM-2017-081. The second named author has been supported by NSF grants DMS-2131243 and DMS--2200762.}
\address{
Department of Mathematics \\
Northeastern University \\
Boston, MA 02115 \\
U.S.A.}
\email{m.yakimov@northeastern.edu}
\date{}
\begin{abstract} The two main approaches to the study of irreducible representations of orders
(via traces and Poisson orders)
have so far been applied in a completely independent fashion. We define and study a 
natural compatibility relation between the two approaches leading to the notion of Poisson trace orders. 
It is proved that all regular and reduced traces are always compatible with any Poisson order structure.  
The modified discriminant ideals of all Poisson trace orders are proved to be Poisson ideals 
and the zero loci of discriminant ideals are shown to be unions of symplectic cores, under natural assumptions 
(maximal orders and Cayley--Hamilton algebras). A base change theorem for Poisson trace orders is proved. 
A broad range of Poisson trace orders are constructed based on the proved theorems: 
quantized universal enveloping algebras, 
quantum Schubert cell algebras and quantum function algebras at roots of unity, 
symplectic reflection algebras, 3 and 4-dimensional Sklyanin algebras, 
Drinfeld doubles of pre-Nichols algebras of diagonal type, and root of unity quantum cluster algebras. 
\end{abstract}
\maketitle
%%%%%%%%%%%%%%%%%%%%   Introduction   %%%%%%%%%%%%%%%%%%%%%%%%%%%%%%%%%%%%%%%%
\sectionnew{Introduction}
\subsection{Traces vs Poisson orders}
There are two main approaches to the study of irreducible representations of orders 
(throughout we will assume that they are algebras over a commutative base ring $\kk$):
\smallskip

(1) The first approach is based on {\em{invariant theory}} and relies on a {\em{trace map}} 
\[
\tr : R \to C, \; \mbox{where $C$ is a central subalgebra}. 
\]
There are several general constructions of traces. If $R$ is a free $C$-module of finite rank, then a $C$-basis of $R$ gives rise to 
a homomorphism from $R$ to the algebra of $n \times n$ matrices over $C$, which, in turn, produces the {\em{regular trace}}
\[
\tr_{\reg} : R \to C.
\]
A more frequently encountered situation is when $R$ is a maximal order,
in which case one has the {\em{reduced trace}}
\[
\tr_{\red} : R \to \ZZ(R),
\]
where $\ZZ(R)$ denotes the center of $R$. 

A {\em{Cayley--Hamilton algebra}} of degree $d \in \Zset_+$ in the sense of Procesi \cite{P} is a $\kk$-algebra $R$ 
with trace $\tr : R \to C$ for a central subalgebra $C$ such that 
\[
\chi_{d,a}(a)=0 \; \;  \mbox{for all} \; \; a\in R \quad \mbox{and} \quad \tr(1)=d,
\]
see \deref{CHring} below. Here $\chi_{d,a}(t) \in C[t]$ is the $d$-th characteristic polynomial of $a \in R$ and 
one requires $1, \ldots, d$ to be invertible in $C$. Every maximal order $R$ of PI degree $d$ 
such that $1, \ldots, d$ are invertible in $R$ is a Cayley--Hamilton algebra of degree $d$ over 
its full center $\ZZ(R)$ with respect to the reduced trace $\tr_{\red}$.  

Using invariant theory, for each $\mm \in \maxspec C$, one constructs an $R/\mm R$-module $V_\mm$ of dimension 
$d$, independent on $\mm$, such that $V_\mm$ is a direct sum (with multiplicities) of all irreducible $R$-modules annihilated by $\mm$. 
For maximal orders $R$ this was done by Braun in \cite{B}; here $d$ equals the PI degree of $R$.
For Cayley--Hamilton algebras of degree $d$ this was done by Procesi in \cite {P}. In this way, 
one has a uniform invariant theory model that captures all irreducible representations of $R$. 
\smallskip

(2) The second approach is based on {\em{Poisson geometry}} and relies on the notion of {\em{Poisson order}}. 
A Poisson order is a triple $(R,C, \partial)$ where $C$ is a central subalgebra of $R$ such that $R$ is a finitely
generated $C$-module and $\partial : C \to \Der_\kk(R)$ is a $\kk$-linear map such that 
\[
\{a, b \}:= \partial_a (b) \in C, \forall a, b \in C  
\]
and $C$ is a Poisson algebra with this operation, see \deref{Pord}. The notion was axiomatized 
in \cite{BG} after the work of De Concini--Kac--Procesi \cite{DKP1} who used it in the framework of big quantum groups 
at roots of unity. Poisson orders are very common and arise whenever we can realize the algebra $R$ as a specialization,
see e.g. \cite[$\S$2.2]{BG} and \cite[Proposition 2.7]{WWY1}.
Their key application is \cite[Theorem 4.2]{BG} stating that the algebras
\[
R/\mm R 
\]
are isomorphic to each other across the symplectic cores of $\maxspec C$ with respect to the Poisson structure $\{.,.\}$.
This drastically simplifies the study of the irreducible representations of $R$ since one only needs to consider 
one finite-dimensional algebra of the form $R/\mm R$ for each symplectic core of $C$. 

For many orders $R$, both approaches (1) and (2) were applied to obtain valuable
information about the irreducible representations of $R$, using invariant theory and Poisson geometry, respectively. 
However, in each case, those applications were carried out independently of each other. 
\subsection{Statements of main results}
The goal of this paper is to unify the two approaches on the basis of the following compatibility condition:
\medskip

\noindent
{\bf{Definition.}} A {\em{Poisson trace order}} is a Poisson order $(R,C, \partial)$ equipped with a trace map $\tr : R \to C$ such that 
\[
\tr \circ \partial_c = \partial_c \circ \tr, \quad \forall c \in C,
\] 
or equivalently, 
\[
\tr( \partial_c(r) ) = \{ c, \tr(r) \}, \quad \forall c \in C, r \in R.
\]

Our first theorem provides a general way of constructing Poisson trace orders based on regular and reduced traces:
\medskip

\noindent
{\bf{Theorem A.}} {\em{Let $(R, C, \partial)$ be a Poisson order.}}
\begin{enumerate}
\item {\em{If $R$ is a free $C$ module, then the quadruple $(R, C, \partial, \tr_{\reg})$ is a Poisson trace oder.}}
\item {\em{If $R$ is a PI algebra and $C = \ZZ(R)$ is normal, then  
$(R, C, \partial, \tr_{\red})$ is a Poisson trace order.}} 
\end{enumerate}
\medskip

\noindent
In fact, we prove a stronger form of the results in both parts of the theorem. Namely:
\begin{enumerate}
\item In the setting of part (1), $\tr_{\reg} \circ \delta = \delta \circ \tr_{\reg}$ for every derivation  
$\delta \in \Der_\kk(R)$ such that $\delta(C) \subseteq C$. 
\item In the setting of part (2), $\tr_{\red} \circ \delta = \delta \circ \tr_{\red}$ for every derivation  
$\delta \in \Der_\kk(R)$ such that $\delta(\ZZ(R)) \subseteq \ZZ(R)$. 
\end{enumerate}
This is proved in Theorems \ref{treg-trace2} and \ref{ttr-comm}. 

The Poisson trace orders constructed from the second part of Theorem A are usually associated to a singular algebra $\ZZ(R)$.
It is advantageous to have a base change theorem which can be used for the construction of new Poisson trace orders from old ones, 
thus leading to Poisson trace orders $(R, C, \partial, \tr)$ with a smooth central subalgebra $C$. This is done in our second theorem:
\medskip

\noindent
{\bf{Theorem B.}} {\em{Assume that
\[
(R, C, \partial, \tr)
\]
is a Poisson trace order with $C$ an integral domain. If $A$ is a Poisson subalgebra of $C$, 
which is normal considered as a commutative algebra, with $C$ finite over $A$, then 
\[
(R,A, \partial|_A, \tr_{C/A} \circ \tr )
\]
is a Poisson trace order, where $\tr_{C/A} : C \to A$ is the trace function of the 
finite (commutative ring) extension $C/A$.}}
\medskip

Using Theorems A and B we construct an extensive list of Poisson trace orders:
\begin{enumerate}
\item All big quantized universal enveloping algebras $\UU_\ep(\g)$ at roots of unity for complex simple Lie algebras $\g$ \cite{DK} 
with respect to both their full centers and their De Concini--Kac--Procesi (smooth) central Hopf subalgebras \cite{DKP1}.
\item All quantum Schubert cell algebras $\UU_\ep^w$ at roots of unity \cite{DKP2}, where $\g$ is an arbitrary 
complex simple Lie algebra and $w$ is a Weyl group element, 
with respect to both their full centers and their De Concini--Kac--Procesi (smooth) central subalgebras \cite{DKP2}.
\item All quantized coordinate rings $F_\ep[G]$ of connected simply connected complex simple algebraic groups $G$ \cite{DL}, 
with respect to both their full centers and their De Concini--Lyubashenko (smooth) central Hopf subalgebras \cite{DL}.
\item All 3 and 4-dimensional Skyanin algebras corresponding to a finite order automorphism of an elliptic curve 
\cite{ATVdB1,ATVdB2,SS} with respect to their full centers \cite{AST,ST}. 
\item The quantum doubles of the bozonizations of all distinguished pre-Nichols algebras \cite{A} of diagonal type 
with finite root systems \cite{H} that belong to a one-parameter family. This large family
contains all contragredient big quantum super-groups at roots of unity. 
The Poisson trace order is 
with respect to their smooth central Hopf subalgebras constructed in \cite{AAY}. 
\item All symplectic reflection algebras \cite{EG} when the parameter $t=0$, with respect 
to their full centers. 
\item All root of unity quantum cluster algebras $\U_\ep(\B)$ \cite{FZ,HLY,MNTY}, 
with respect to their full centers and with respect to the central subalgebras 
associated to their quantum Frobenuis maps. 
\end{enumerate}
We refer the reader to Sections \ref{5} and \ref{6} for the statements of these results and the precise (mild) assumptions 
in each case. 

Our last theorem addresses the Poisson properties of the discriminant and modified discriminant ideals of Poisson trace orders.
\medskip

\noindent
{\bf{Theorem C.}} The following hold for every  Poisson trace order $(R,C, \partial, \tr)$:
\begin{enumerate}
\item All modified discriminant ideals $MD_\ell(R/C, \tr)$ of $R$ are Poisson ideals of $C$ 
with respect to the underlying Poisson structure on $C$. 
\item If the base ring $\kk$ is a field, $C$ is a finitely generated $\kk$-algebra and 
\begin{enumerate}
\item $R$ is a maximal order, $C = \ZZ(R)$ and $tr : R \to C$ is the reduced trace or 
\item $(R,C, \tr)$ is a Cayley--Hamilton algebra of degree $d$ such that $\charr \kk \notin [1,d]$, 
\end{enumerate}
then the zero sets of all discriminant ideals $D_\ell(R/C, \tr)$ are unions of symplectic cores of $\maxspec C$.  
\end{enumerate}
\medskip
\subsection{Organization of the paper and notation}
The paper is organized as follows. Section \ref{2} contains background material on algebras with traces and Poisson orders, 
defines Poisson trace orders, proves Theorem C and the compatibility property of the regular trace in part (1) of Theorem A. 
Section \ref{3} contains results on the reduced trace map of PI algebras and proves the compatibility property in part (2) of Theorem B. 
Section \ref{4} proves results on Poisson orders coming from extensions of Poisson algebras and the base change 
Theorem B for Poisson orders. Section \ref{5} contains the construction of a broad range of families of Poisson orders with 
respect to the full centers of PI algebras. Section \ref{6} contains the construction of many families of Poisson orders of different nature with 
respect to proper central subalgebras.
\medskip

\noindent
{\bf{Notation:}} {\em{Throughout the paper $\kk$ will denote a base commutative ring. By an algebra we will mean a $\kk$-algebra,
and by a Poisson algebra, a Poisson algebra over $\kk$}}, i.e., the Poisson bracket of the latter is $\kk$-linear.
The center of a $\kk$-algebra $R$ will be denoted by $\ZZ(R)$.
%%%%%
\sectionnew{Poisson trace orders and their discriminant ideals}
\label{2}
In this section we define Poisson trace orders and prove a theorem that their modified discriminant ideals are always Poisson ideals. 
From this we deduce that the zero sets of their discriminant ideals are unions of symplectic cores. 
Furthermore we prove that the regular trace is compatible with any Poisson order structure.

\subsection{Definitions}
We first recall the definitions of Poisson order and trace map for a noncommutative algebra.
Throughout the section $C \subseteq R$ will denote a central $\kk$-subalgebra of $R$ 
over which $R$ is a finitely generated module.
If $R$ is an affine $\kk$-algebra and $\kk$ is a Noetherian commutative ring, 
then by the Artin--Tate Lemma \cite[Lemma 13.9.10]{McR} $C$ 
is an affine $\kk$-algebra, and so, $C$ and $R$ are Noetherian algebras.

\bde{Pord} \cite[Definition 2.1]{BG} The pair $(R, C)$ is called a {\em{Poisson order}} if 
\begin{enumerate}
\item $R$ is a finitely generated $C$-module;
\item $C$ is equipped with a Poisson algebra structure $\{.,.\}$;
\item there exists a map $\partial : C \longrightarrow \Der_\kk(R)$, $c \mapsto \partial_c$ such that $\partial_c$ is an extension of the 
Hamiltonian derivation $\{ c, - \}$ of $C$. 
\end{enumerate}
\ede

We will denote Poisson orders as triples $(R,C, \partial)$ because the Poisson bracket $\{.,.\}$ can be recovered from the map 
$\partial : C \longrightarrow \Der_k(R)$ by condition (iii).

\bde{trace} A {\em{trace map}} from $R$ to $C$ is a map
$\tr \colon R \to C $ that has the following properties:
\begin{enumerate}
\item ($C$-linearity) $\tr(zx) = z \tr(x)$ for $x \in R$, $z \in C$;
\item (cyclicity) $\tr( xy) = \tr(yx)$ for $x, y \in R$.
\end{enumerate}
\ede
\bre{nontriviality} For some purposes, one needs to require that $\tr \neq 0$, or more specifically, that $\tr(1) \in \kk^*$, see e.g. \cite[Definition 2.1]{BY}.
This will not be needed in this paper. 
\ere

Our main definition combines the above two notions:

\bde{P-tr-ord} A {\em{Poisson trace order}} is a Poisson order $(R,C, \partial)$ equipped with a trace map $\tr : R \to C$ such that 
\[
\tr \circ \partial_c = \partial_c \circ \tr, \quad \forall c \in C.
\] 
\ede  
More explicitly, the compatibility condition states that
\[
\tr( \partial_c(r) ) = \{ c, \tr(r) \}, \quad \forall c \in C, r \in R.
\]
\subsection{Poisson properties of modified discriminant ideals}
Our first result is that the modified discriminant ideals of a Poisson trace order behave well with respect to the underlying Poisson structure. 
First, recall the following:

\bde{discr} \cite[Definition 1.2(2)]{CPWZ} For an algebra $R$ with trace $\tr : R \to C$ and $\ell \in \Zset_+$, 
the \emph{modified $\ell$-discriminant ideal} of $R$ with respect to the trace map $\tr$ 
is the ideal $MD_\ell(R/C, \tr)$ of $C$ generated by
\[
d_\ell( \{ r_1, \ldots, r_\ell \}, \{s_1, \ldots, s_\ell \} : \tr ) := \det \big( \tr(r_i s_j)_{i,j =1}^\ell \big)
\] 
for all $\ell$-tuples 
\[
\{ r_1, \ldots, r_\ell \}, \{s_1, \ldots, s_\ell \}  \in R^\ell.
\]
\ede

\bth{Ptrace-Pid} Let $(R,C, \partial, \tr)$ be a Poisson trace order and fix $\ell \in \Zset_+$. Then the modified discriminant ideals $MD_\ell(R/C, \tr)$  are 
Poisson ideals of $C$ with respect to the underlying Poisson structure on $C$. 
\eth
The theorem follows from of the next more general result.
\bpr{der-det} Assume that $(R,C, \tr)$ is an algebra with trace and $\delta$ is a $\kk$-derivation of $R$ that preserves $C$ and commutes with $\tr$: 
\[
\delta \big(  \tr (r)  \big)= \tr \big( \delta (r) \big), \quad \forall r \in R. 
\]
Then for all $\ell$-tuples 
\[
\{ r_1, \ldots, r_\ell \}, \{s_1, \ldots, s_\ell \}  \in R^\ell,
\]
we have 
\begin{align*}
\delta d_\ell( \{ r_1, \ldots, r_\ell \}, \{s_1, \ldots, s_\ell \} : \tr ) 
&= \sum_{k=1}^\ell d_\ell( \{ r_1, \ldots, \delta r_k, \ldots r_\ell \}, \{s_1, \ldots, s_\ell \} : \tr ) \\
&+ \sum_{k=1}^\ell d_\ell( \{ r_1, \ldots, r_\ell \}, \{s_1, \ldots, \delta s_k, \ldots s_\ell \} : \tr ). 
\end{align*}
\epr
\begin{proof} We have
\begin{align*}
&\delta d_\ell( \{ r_1, \ldots, r_\ell \}, \{s_1, \ldots, s_\ell \} : \tr ) 
= \sum_{\sigma \in S_\ell} (-1)^{\sign \sigma} \delta \Big( \prod_{k=1}^\ell \tr(r_k s_{\sigma(k)} ) \Big)
\\
=& \sum_{k=1}^\ell \sum_{\sigma \in S_\ell} (-1)^{\sign \sigma} \tr(r_1 s_{\sigma(1)} ) \ldots \tr(\delta (r_k) s_{\sigma(k)} ) \ldots \tr(r_\ell s_{\sigma(\ell)} ) 
\\
+& \sum_{k=1}^\ell \sum_{\sigma \in S_\ell} (-1)^{\sign \sigma} \tr(r_1 s_{\sigma(1)} ) \ldots \tr( r_k \delta (s_{\sigma(k)}) ) \ldots \tr(r_\ell s_{\sigma(\ell)} ) 
\\
=&\sum_{k=1}^\ell d_\ell( \{ r_1, \ldots, \delta r_k, \ldots r_\ell \}, \{s_1, \ldots, s_\ell \} : \tr ) \\
+&\sum_{k=1}^\ell d_\ell( \{ r_1, \ldots, r_\ell \}, \{s_1, \ldots, \delta s_k, \ldots s_\ell \} : \tr ). 
\end{align*}
\end{proof} 
\subsection{Poisson properties of discriminant ideals}
The modified discriminant ideals of an algebra with trace defined in \cite{CPWZ} 
are newer versions of the much older (and now classical) notion of discriminant ideals:

\bde{discr0} \cite[p. 126]{Re} For an algebra $R$ with trace $\tr : R \to C$ and $\ell \in \Zset_+$, 
the \emph{$\ell$-discriminant ideal} of $R$ with respect to $\tr$ 
is the ideal $D_\ell(R/C, \tr)$ of $C$ generated by
\[
d_\ell( \{ r_1, \ldots, r_\ell \}, \{r_1, \ldots, r_\ell \} : \tr ) := \det \big( \tr(r_i r_j)_{i,j =1}^\ell \big)
\] 
for all $\ell$-tuples 
\[
\{ r_1, \ldots, r_\ell \} \in R^\ell.
\]
\ede

While we cannot prove or disprove that the discriminant ideals of Poisson trace orders are Poisson ideals, 
we can show that their zero loci are unions of symplectic cores under natural assumptions that are satisfied in 
wide generality. 

Symplectic cores are algebraic versions of symplectic leaves that can be defined for the 
maximal spectra of finitely generated Poisson algebras in any characteristic unlike symplectic leaves \cite[Sect. 5]{CW}, which are 
defined via transcendental methods when base field is $\Rset$ or $\Cset$. Let $(C, \{.,.\})$ be a Poisson algebra over a field $\kk$,
which is finitely generated as a commutative algebra.
For every ideal $I$ of $C$, there exists a unique maximal Poisson ideal $\PP(I)$ contained in $I$; $\PP(I)$ is Poisson prime when $I$ is prime 
\cite[Lemma~6.2]{G}. Following \cite[Section 3.2]{BG}, the {\em{symplectic core}} of $\mm \in \maxspec C$ is the subset of 
$\maxspec C$ given by
\[
\CC(\mm) := \{ \n \in \maxspec C \mid \PP(\n) = \PP(\mm) \}.
\]
Symplectic cores define a partition of $\maxspec C$ by locally closed subsets.

For our results on discriminant ideals, we need to recall Procesi's notion of Cayley--Hamilton algebras of degree $d$ \cite {P}, where 
$d$ is a positive integer.
Consider an algebra with trace $(R, C, \tr)$ over a field $\kk$ of characteristic $0$ or $>d$.
For $1 \leq k \leq d$, denote by $\sigma_k$ the $k$-th elementary symmetric function in the indeterminates 
$\lambda_1, \lambda_2, \dots, \lambda_d$ and by $\psi_k:=\lambda_1^k+\lambda_2^k+\cdots +\lambda_d^k$ the $k$-th Newton power sum function.
As is well known, that there exists a unique set of polynomials 
\[
p_i(x_1, x_2, \dots, x_k)\in \mathbb{Z}[(k!)^{-1}][x_1, x_2, \dots, x_k]
\] 
such that
\[
\sigma_k=p_k(\psi_1, \psi_2, \dots, \psi_k), \quad \forall \; 1\leq k \leq d.
\]
The $d$-th characteristic polynomial $\chi_{d,a}(t) \in C[t]$ 
of an element $a \in R$ is defined by
\[
\chi_{d, a}(t):= t^d-c_1(a)t^{d-1}+\cdots +(-1)^d c_d(a),
\]
where $c_k(a):=p_k\big(\tr(a), \tr(a^2), \dots, \tr(a^k)\big)$.

\bde{CHring} \cite{P} A Cayley--Hamilton algebra of degree $d \in \Zset_+$ is a $\kk$-algebra with trace $(R, C, \tr)$
such that
\begin{enumerate}
\item $\chi_{d,a}(a)=0$ for all $a\in R$ and
\item $\tr(1)=d$.
\end{enumerate}
\ede

\bco{Ptrace-Pid2} Let $(R,C, \partial, \tr)$ be a Poisson trace order over a field $\kk$ such that 
$C$ is a finitely generated $\kk$-algebra and either
\begin{enumerate}
\item $R$ is a maximal order, $C = \ZZ(R)$ and $tr : R \to C$ is the reduced trace {\em{(}}see Sect. \ref{3}{\em{)}} or 
\item $(R,C, \tr)$ is a Cayley--Hamilton algebra of degree $d$ such that $\charr \kk \notin [1,d]$. 
\end{enumerate}
Then, for all $\ell \in \Zset_+$, the zero set of the discriminant ideals $D_\ell(R/C, \tr)$ is a union 
of symplectic cores of $\maxspec C$.  
\eco
\begin{proof} Main Theorem (a)-(b) and Theorem 4.1(b) in \cite{BY} imply that under the assumption (1) or (2), 
the zero set of the discriminant ideal $D_\ell(R/C, \tr)$ coincides with that of the 
modified discriminant ideal $MD_\ell(R/C, \tr)$. The latter is a Poisson ideal 
of $(C, \{.,.\})$ by \thref{Ptrace-Pid}, and hence its zero locus is a union of symplectic 
cores of  $\maxspec C$. 
\end{proof}
The conclusion of \coref{Ptrace-Pid2} is valid under the more general assumption that $\tr : R \to C$ is a 
{\em{representation theoretic trace}} in the sense of \cite[Definition 2.1]{BY}; this follows by applying 
\cite[Main Theorem (a)]{BY}. 
\subsection{Compatibility properties of the regular trace map}
We finish the section with a general construction of Poisson trace orders from regular trace maps.
Assume that $R$ is a free module of finite rank over the central subalgebra $C$. The {\em{regular trace}} is the composition
\begin{equation}
\label{reg-tr}
\tr_{\reg} : R \to \End_C(R) \cong M_n(C) \stackrel{\tr}{\lra} C.
\end{equation}
The first map is given by left multiplication, the second one is the isomorphism obtained by choosing a  
$C$-basis $\{v_1, \ldots v_n \}$ of $R$, and the third one is the matrix trace. The map is obviously independent on 
the choice of $C$-basis of $R$.

\bth{reg-trace1} Assume that $R$ is an algebra and $C$ is a central subalgebra of $R$
such that $R$ is a finite rank free $C$-module. Then every structure of Poisson order $(R, C, \partial)$ is a Poisson trace order with 
respect to the regular trace $\tr_{\reg} : R \to C$.
\eth

In other words, the required compatibility between a Poisson order map $\partial : C \to \Der_\kk(R)$ and a trace map $\tr : R \to C$ 
is always satisfied for regular trace maps. 

The theorem follows immediately from the following theorem.

\bth{reg-trace2} If, in the setting of \thref{reg-trace1}, 
\[
\delta \in \Der_\kk(R) \;\; \mbox{is such that} \; \; \delta(C) \subseteq C, 
\]
then 
\[
\tr_{\reg} \circ \delta = \delta \circ \tr_{\reg}.  
\]
\eth
A proof of this fact was given in \cite[Proposition 2.2]{NTY} using deformation theory. We next give a direct proof.
\begin{proof} Fix a $C$-basis $\{v_1, \ldots v_n \}$ of $R$. Let $r \in R$.
Denote
\[
r v_i = \sum_j b_{ij} v_j
\] 
and
\[
\delta(v_j) = \sum_k c_{jk} v_k
\]
for some $b_{ij}, c_{jk} \in C$. We have
\begin{align*}
\delta (r) v_i &= \delta (r v_i) - r \delta (v_i)
\\
&= \sum_j \delta (b_{ij}) v_j  + \sum_j b_{ij} \delta (v_j) - r \delta (v_i)
\\
&= \sum_j \delta (b_{ij}) v_j  + \sum_j b_{ij} \delta (v_j) - \sum_k c_{ik} r v_k 
\\
&= \sum_j \delta (b_{ij}) v_j  + \sum_{j,k} b_{ij} c_{jk} v_k - \sum_{k,\ell} c_{ik} b_{k\ell}.
\end{align*}
Therefore, 
\begin{align*}
\tr_{\reg} \big( \delta(r) v_i \big) &= \sum_i \delta(b_{ii}) + \sum_{i,j} b_{ij} c_{ji} - \sum_{k,i} c_{ik} b_{ki} = 
\\
&=  \sum_i \delta(b_{ii}) =  \delta \big( \tr_{\reg} (r) \big).
\end{align*}
\end{proof}
%%%%%%%%
\sectionnew{Poisson trace orders for maximal orders}
\label{3}
In this section we prove that the reduced trace on a prime PI algebra whose center is normal is compatible with any Poisson order structure.

\subsection{Reduced traces of PI algebras}
We first recall the definition of reduced trace for PI algebras. Let $R$ be a prime PI ring, see \cite[Definition I.13.1]{BG-book}. 
Denote by $Q$ the fraction field of its center $\ZZ(R)$. The ring of fractions of $R$ is isomorphic to the tensor product
\[
R \otimes _{\ZZ(R)} Q.
\]
By Posner's theorem, \cite[$\S$I.13.3]{BG-book}, it is a 
central simple $Q$-algebra. Denote by $n$ the PI degree of $R$; that is, $\dim_Q (R \otimes_{\ZZ(R)} Q) = n^2$, \cite[$\S$I.13.3]{BG-book}.
Since $R$ is a finitely generated torsion-free $\ZZ(R)$-module, it is, in the terminology of \cite[$\S$8]{Re}, a
\emph{$\ZZ(R)$-order} in $R \otimes _{\ZZ(R)} Q$. Denote the embedding
\[
\iota : R \hra R \otimes _{\ZZ(R)} Q.
\]
There exists a finite field extension $F$ of $Q$ which splits
$R \otimes _{\ZZ(R)} Q$; that is,
\[
R \otimes_{\ZZ(R)} F \cong M_n(F),
\]
see  \cite[\S 7b]{Re}. The field embedding $\mu : Q \hra F$ gives rise 
to the ring homomorphism
\begin{equation}\label{added}
\phi : = \id_R \otimes \mu : R \otimes_{\ZZ(R)} Q \to R \otimes_{\ZZ(R)} F.
\end{equation}

Consider the composition
\[
R \stackrel{\iota}{\hookrightarrow} R \otimes_{\ZZ(R)} Q \stackrel{\phi}{\lra} R \otimes_{\ZZ(R)} F \cong M_n(F) \stackrel{\tr}{\lra} F,
\]
where the last map is the matrix trace. This composition takes values in $Q \subset F$, 
\[
\tr_{\red} : R \to Q
\]
and is called the {\em{reduced trace map}} of $R$. If $\ZZ(R)$ is normal, the reduced trace map takes values in $\ZZ(R) \subset Q$, see \cite[Theorem 10.1]{Re}:
\[
\tr_{\red} : R \to \ZZ(R).
\]

Since $R$ is prime, its center $\ZZ(R)$ is an integral domain. If $\charr R \nmid n$, then
\[
\tr_{\red}(1) = n \neq 0,
\]
cf. \reref{nontriviality}. 
\subsection{Compatibility properties of the reduced trace map} Our next result is a general compatibility theorem between reduced traces and Poisson orders. 
\bth{stan-red}
Let $R$ be a $\kk$-algebra over a commutative ring $\kk$ such that 
\begin{enumerate}
\item $R$ is a prime algebra and
\item the center $\ZZ(R)$ is normal.
\end{enumerate}
Then every structure of Poisson order $(R, \ZZ(R), \partial)$ is a Poisson trace order with respect to the reduced trace $\tr_{\red} : R \to \ZZ(R)$. 
\eth

The point of the theorem is that the needed compatibility between the Poisson order map $\partial : \ZZ(R) \to \Der_\kk(R)$ and 
the trace map $\tr : R \to \ZZ(R)$ is automatically satisfied for a reduced trace map. 

We have an important special case of \thref{stan-red}: if $R$ is a maximal order, then condition (2) is satisfied, by \cite[Proposition 5.1.10(b)]{McR}.

The theorem follows from the next more general result on the commutation of the reduced trace with any derivation 
of the algebra $R$ that preserves (but does not necessarily fix) its center $\ZZ(R)$. 

\bth{tr-comm} Assume the setting of \thref{stan-red}. Let $\delta \in \Der(R)$ be such that $\delta(\ZZ(R)) \subseteq \ZZ(R)$.
Then 
\[
\tr_{\red} \circ \delta = \delta \circ \tr_{\red}
\]
\eth
\noindent
\begin{proof}
Recall from the beginning of the section that $Q$ denotes the field of fractions of $\ZZ(R)$. The derivation $\delta$ of $R$ uniquely extends 
to a derivation of the quotient ring $R \otimes_{\ZZ(R)} Q$,
\[
\delta' \in \Der( R \otimes_{\ZZ(R)} Q);
\]
that is,
\begin{equation}
\label{delta-iota}
\delta' \circ \iota = \iota \circ \delta.
\end{equation}
Furthermore, since $\delta( \ZZ(R)) \subseteq \ZZ(R)$, we have $\delta'(Q) \subseteq Q$.
Because the field extension $F/Q$ can be chosen to be finite and separable \cite[Theorem 7.10]{Re}, 
$\delta'$ extends to a derivation $\delta''_0$ of $F$ \cite[Theorem VIII.5.1]{L}, which in turn extends to a derivation $\delta''$ of 
\[
R \otimes_{\ZZ(R)} F \cong M_n(F) \quad \mbox{by} \quad \delta''(r \otimes f): = \delta(r) \otimes f + r \otimes \delta''_0(f), \; \; \forall r \in R, f \in F.
\]
This derivation satisfies $\delta''(F) \subseteq F$ and, keeping the notation (\ref{added}),
\begin{equation}
\label{delta-phi}
\delta'' \circ \phi = \phi \circ \delta'.
\end{equation}
Obviously the map
\[
\ol{\delta}'' : M_n(F) \to M_n(F), \; \;  \mbox{given by} \; \; 
\ol{\delta}''(f E_{ij}) := \delta''(f) E_{ij}, \; \; \forall 1, \leq i,j \leq n,\, \forall f \in F,
\]
is a derivation of $M_n(F)$ and $\delta'' - \ol{\delta}''$ is an $F$-linear derivation of $M_n(F)$:
\[
\delta'' - \ol{\delta}'' \in \Der_F( M_n(F)).
\]
It is a corollary of the Skolem--Noether theorem  that every $F$-linear derivation of a central simple $F$-algebra is inner; see \cite[Theorem 7.21]{Re}, also \cite[Proposition 1]{Am}. Thus, $\delta'' - \ol{\delta}''$ is an inner derivation of $M_n(F)$. So, there exists $X \in M_n(F)$ such that 
\[
\delta'' = \ol{\delta}'' + [X, -]. 
\]
Therefore,
\[
\delta'' (f E_{ij}) := \delta''(f) E_{ij} + [X, f E_{ij}], \quad \forall 1 \leq i,j \leq n, f \in F.
\] 
Let $\delta_{ij}$ denote the Kronecker delta. For the matrix trace $\tr: M_n(F) \to F$, we have 
\[
\tr \big( \delta'' (f E_{ij}) \big) = \tr \big( \delta''(f) E_{ij} \big)  + \tr \big( [X, f E_{ij}] \big) = \delta_{ij} \delta''(f)
\]
and 
\[
\delta'' \big(  \tr (f E_{ij}) \big) = \delta'' ( \delta_{ij} f ) = \delta_{ij} \delta''(f)
\]
for all $1 \leq i,j \leq n$, $f \in F$. Therefore, 
\begin{equation}
\label{tr-delta}
\tr \circ \delta'' = \delta'' \circ \tr.
\end{equation}
Now $\ZZ(R)$ is a subring of $Q$ via $\iota$ and $Q$ is a subfield of $F$ via $\phi$, so $\iota|_{\ZZ(R)} = \id_{\ZZ(R)}$, 
$\phi|_Q = \id_Q$. It follows from \eqref{delta-iota}--\eqref{delta-phi} that for all $r \in R$, 
\[
\delta'' \circ \phi \circ \iota (r) =  \phi \circ \delta' \circ \iota (r) =   \phi \circ \iota \circ \delta (r).
\]
In particular, 
\[
\delta''(z) = \delta (z), \quad \forall z \in \ZZ(R). 
\]
Combining the last two identities, eq. \eqref{tr-delta}, and the definition of the reduced trace $\tr_{\red} : R \to \ZZ(R)$, 
gives that for all $r \in R$,
\begin{align*}
\tr_{\red} \big( \delta (r) \big)&= \tr \circ \phi \circ \iota \circ \delta(r) 
\\
&=\tr \circ \delta'' \circ \phi \circ \iota  (r)
\\
&= \delta'' \circ \tr \circ \phi \circ \iota (r)
\\
&= \delta'' \big( \tr_{\red} (r) \big) = \delta \big( \tr_{\red}(r) \big).
\end{align*} 
Therefore, $\tr_{\red} \circ \delta = \delta \circ \tr_{\red}$. 
\end{proof}
%%%%%%%%
\sectionnew{Base change for Poisson trace orders}
\label{4}
In this section we prove a base change theorem on Poisson trace orders. 

\subsection{A setting for base change for Poisson trace orders}
In applications to the representation theory of PI algebras $R$, it is advantageous to work with Poisson orders $(R,C, \partial)$ 
for nonsingular algebras $C$ because one has to deal with the Poisson geometry of its spectrum 
\cite{DKP1,DKP2,BG,WWY1,WWY2}.
The Poisson trace orders $(R, Z(R), \partial, \tr)$ arising from \thref{stan-red} often involve singular centers $Z(R)$.
This raises the problem of analysing how to relate a Poisson trace order $(R, C, \partial, \tr)$ with an order $(R,A, \partial|_A, \tr')$, 
where $A$ is a Poisson subalgebra of $C$ such that $C$ is finite over $A$ (i.e., $C$ is a finitely generated $A$-module). That is, symbolically, we examine the change
\begin{equation}
\label{base-change}
(R, C, \partial, \tr) \rightsquigarrow (R, A, \partial|_A, \tr').
\end{equation}
 In this situation $(R, A, \partial|_A)$ is necessarily a Poisson order. If $\tr_{C/A} : C \to A$ is a trace map of commutative algebras, then we have the composition 
trace map 
\[
\tr' := \tr_{C/A}\circ \tr : R \to A,
\]
which is the map that is used in the base change in \eqref{base-change}. The question addressed in this section is: 
\medskip

{\em{Under what conditions is $(R, A, \partial|_A, \tr')$ again a Poisson trace order?}}

\subsection{Poisson algebra extensions}
To deal with base changes, we first investigate Poisson trace orders arising from Poisson algebra extensions. Let $A$ be a Poisson subalgebra of 
a Poisson algebra $(C,\{.,.\})$ such that $C$ is finite over $A$. Then we have an induced Poisson order structure on the pair $(C,A)$ with 
the map $\partial^{\res}$ given by
\begin{equation}
\label{restrict-part}
\partial^{\res}_a := \{ a , . \}, \quad \forall a \in A.
\end{equation} 
Assume that $C$ is an integral domain, and denote by $\FF(C)$ and $\FF(A)$ the fields of fractions of $C$ and $A$, respectively. Then 
$\FF(C)/\FF(A)$ is a finite field extension and we have the regular trace function
\begin{equation}
\label{tr-FF}
\tr_{\FF(C)/\FF(A)} : \FF(C) \to \FF(A).
\end{equation}
If $A$ is normal as a commutative algebra, then its restriction to $C$ takes values in $A$
\begin{equation}
\label{tr-res}
\tr_{C/A} := \tr_{\FF(C)/\FF(A)}|_{C} : C \to A.
\end{equation}
This is so because for all $c \in C$, $\tr_{\FF(C)/\FF(A)}(a)$ is integral over $A$ and belongs to $\FF(A)$, so by the normality assumption on $A$, it 
belongs to $A$. 
\bpr{comm-alg-ext} Assume that $(C, \{.,.\})$ is a Poisson algebra and an integral domain. 
Let $A$ be a Poisson subalgebra, which is normal considered as a commutative algebra, 
with $C$ finite over $A$. Then 
\[
(C,A, \partial^{\res}, \tr_{C/A})
\]
is a Poisson trace order.
\epr
\begin{proof} The Poisson structure $\{.,.\}$ on $C$ uniquely extends to a Poisson (field) structure on $\FF(C)$ and $\FF(A)$ is a 
Poisson subfield of $\FF(C)$. By abuse of notation, this extension will be denoted by the same symbol. We obtain a Poisson order structure
\[
(\FF(C), \FF(A), \ol{\partial}^{\res})
\]
by setting
\[
\ol{\partial}^{\res}_f := \{f,.\}, \quad \forall f \in \FF(A).
\]
Clearly, we have 
\begin{equation}
\label{part-rest}
\partial^{\res} = \ol{\partial}^{\res}|_A.
\end{equation}
The trace function \eqref{tr-FF} is the regular trace of the pair $(\FF(C), \FF(A))$ given by \eqref{reg-tr}. 
The standard definition of trace function of the field extension $\FF(C)/\FF(A)$ is given by a sum over the field embeddings 
of $\FF(C)$ into the algebraic closure of $\FF(A)$, \cite[Chap. VI, Sec. 5]{L}; and this is equivalent to the regular trace by \cite[Proposition VI.5.6]{L}. 
Now we can apply \thref{reg-trace1} to obtain that 
\[
\big( \FF(C), \FF(A), \ol{\partial}^{\res}, \tr_{\FF(C)/\FF(C)} \big)
\]
is a Poisson trace order. This shows that
\[
\tr_{\FF(C)/\FF(A)} \big( \ol{\partial}^{\res}_f(g) \big) = \{ f, \tr_{\FF(C)/\FF(A)}(g) \}, \quad \forall f \in \FF(A), g \in \FF(C).
\]
Applying this for $f \in A$ and $g \in C$, and combining it with \eqref{tr-res} and \eqref{part-rest}, 
gives
\[
\tr_{C/A} \big( \partial^{\res}_a (c) \big) = \{ a, \tr_{C/A}(c) \}, \quad \forall a \in A, c \in C.
\]
\end{proof}
\subsection{Statement of base change theorem} The following is our base change theorem for Poisson trace orders. 
\bth{base-change} Let 
\[
(R, C, \partial, \tr)
\]
be a Poisson trace order with $C$ an integral domain. Assume that $A$ is a Poisson subalgebra of $C$, 
which is normal considered as a commutative algebra, with $C$ finite over $A$. Then 
\[
(R,A, \partial|_A, \tr_{C/A} \circ \tr )
\]
is a Poisson trace order.
\eth
\begin{proof} The first assumption gives that
\[
\tr( \partial_a(r) ) = \{ a, \tr(r) \}, \quad \forall a \in A, r \in R.
\]
Applying \prref{comm-alg-ext} and using the definition \eqref{restrict-part} of $\partial^{\res}$ gives
\[
\tr_{C/A} \{ a , c \} = \{ a, \tr_{C/A}(c) \},  \quad \forall a \in A, c \in C.
\]
Therefore, 
\[
\tr_{C/A} \big( \tr( \partial_a(r) ) \big) = \tr_{C/A} \big(  \{ a, \tr(r) \} \big) = \{ a, \tr_{C/A}(\tr(r)) \}
\]
for all $a \in A$, $r \in R$,
which completes the proof of the theorem.
\end{proof}
\sectionnew{Construction of Poisson trace orders: full centers}
\label{5}
In this section we use \thref{stan-red} to construct large classes of Poisson trace orders 
with respect to their full centers from the areas of quantum groups,
noncommutative projective algebraic geometry, cluster algebras and Lie theory (symplectic reflection algebras).

\subsection{Specializations plus maximal orders}
\label{5.1}
\thref{stan-red} implies that we have the following very general construction of Poisson trace orders:
\bpr{spec-red}
Assume that $R$ is a prime affine $\kk$-algebra which is finitely generated over its center. 
\begin{enumerate}
\item If $R$ is obtained as a specialization as in \cite[$\S$2.2]{BG} or, more generally, a higher order specialization as in \cite[Definitions 2.3 and 2.5]{WWY1}, 
then by \cite[Proposition 2.7]{WWY1}, $(R, \ZZ(R))$ has a canonical structure of Poisson order $\partial_{\spp}$.
\item If $\ZZ(R)$ is normal, then we can consider the reduced trace map $\tr_{\red} : R \to \ZZ(R)$.
\end{enumerate}
Under the assumptions (1)-(2), $(R, \ZZ(R), \partial_{\spp}, \tr_{\red})$ is a Poisson trace order.
\epr
Condition (2) is satisfied if $R$ is a maximal order \cite[Proposition 5.1.10(b)(i)]{McR}.

\subsection{Quantized enveloping algebras, quantum functions algebras and quantum Schubert cell algebras at roots of unity}
\label{5.2}
Let $\g$ be a complex simple Lie algebra, $G$ be the corresponding connected simply connected complex simple algebraic group, and $w$ be 
an element of the Weyl group of $\g$. Let $\ep \in \Cset$ be a primitive root of unity of odd order $\ell$ such that $3 \nmid \ell$ if $\g$ is of type $G_2$. 
We have the following three affine $\Qset(\ep)$-algebras:
\begin{enumerate}
\item The big quantized universal enveloping algebra of $\g$ at root of unity, denoted by $\UU_\ep(\g)$ and constructed in \cite{DK}.
\item The quantized coordinate ring of $G$ at root of unity, denoted $F_\ep[G]$ and constructed in \cite{DL}.
\item The quantum Schubert cell algebra at root of unity associated to $\g$ and $w$, denoted by $\UU_\ep^w$ and constructed in \cite{DKP2}. 
\end{enumerate} 

The three families of algebras are obtained by specialization as shown in \cite{DKP1}, \cite{DL} and \cite{DKP2}, respectively.
This gives rise to the Poisson order structures
\begin{equation}
\label{PO}
(\UU_\ep(\g), \ZZ(\UU_\ep(\g)), \partial_{\spp}), \quad 
(F_\ep[G], \ZZ(F_\ep[G]), \partial_{\spp}), \quad
(\UU_\ep^w, \ZZ(\UU_\ep^w), \partial_{\spp}).
\end{equation}
The three Poisson orders were shown to be nontrivial, and the Poisson structures on their centers 
were described in terms of Poisson--Lie groups and Poisson structures on Schubert cells in
\cite{DK,DKP1,DKP2,DL}.   

The algebras in (1)--(3) are also maximal orders by \cite[Theorem 1.8]{DK}, \cite[Theorem 7.4]{DL} and \cite[Theorem 1.5]{DKP2}, respectively.
Thus, we have a reduced trace map $\tr_{\red} : R \to \ZZ(R)$ for each algebra $R$ in the three families. 

\prref{spec-red} now implies the following:
\bpr{1} For all complex simple Lie algebras $\g$, connected simply connected complex simple algebraic groups $G$,
Weyl group elements $w$ and primitive roots of unity $\ep \in \Cset$ of odd order $\ell$ such that $3 \nmid \ell$ if $\g$ is of type $G_2$,
the quadruples
\begin{multline*}
(\UU_\ep(\g), \ZZ(\UU_\ep(\g)), \partial_{\spp}, \tr_{\red}), \quad 
(F_\ep[G], \ZZ(F_\ep[G]), \partial_{\spp}, \tr_{\red}) \\
\mbox{and} \quad (\UU_\ep^w, \ZZ(\UU_\ep^w), \partial_{\spp}, \tr_{\red}).
\end{multline*}
are Poisson trace orders:
\epr
\subsection{Noncommutative projective spaces}
\label{5.3}
Sklyanin algebras are (affine) quadratic algebras that play a major role in Skyanin's work on quantum integrable systems \cite{Sk}, 
the Artin--Schelter--Tate--van den Bergh classification of noncommutative projective spaces \cite{AS,ATVdB1,ATVdB2}, and the Feigin--Odeskii 
study of elliptic algebras \cite{FO}. 
%The explicit definitions of these algebras in terms of generators and relations can be found for example 
%in \cite[Definition 2.12]{WWY1} and \cite[Definition 2.1]{WWY2}.

The 3-dimensional Sklyanin algebras are the $\Cset$-algebras in three generators $x, y, z$ subject to relations
\[
ayz+bzy+cx^2 ~=~ azx+bxz+cy^2 ~=~ axy+byx+cz^2 ~=~ 0
\]
for parameters $a, b, c \in \Cset$ such that 
\[
abc \neq 0 \; \; \mbox{and} \; \; (3abc)^3 \neq (a^3+b^3+c^3)^3. 
\]

The 4-dimensional Sklyanin algebra
are the $\Cset$-algebras in four generators 
$x_0$, $x_1$, $x_2$, $x_3$, subject to the following relations  
\[
\begin{array}{llll}
x_0x_1-x_1x_0 ~=~  \alpha(x_2x_3+x_3x_2), &&& x_0x_1+x_1x_0 ~=~ x_2x_3-x_3x_2,\\
x_0x_2-x_2x_0 ~=~  \beta(x_3x_1+x_1x_3), &&& x_0x_2+x_2x_0 ~=~  x_3x_1-x_1x_3,\\
x_0x_3-x_3x_0 ~=~  \gamma(x_1x_2+x_2x_1), &&& x_0x_3+x_3x_0 ~=~  x_1x_2-x_2x_1,\\
\end{array}
\]
for parameters $\alpha, \beta, \gamma \in \Cset$ such that
\[
\alpha + \beta + \gamma + \alpha \beta \gamma = 0  \; \; \mbox{and} \; \; 
(\alpha, \beta), (\beta, \gamma), (\gamma, \alpha) \neq (-1, 1). 
\]

The following properties were proved in \cite{AS,ATVdB1,ATVdB2,SS} for each Sklyanin algebra $S$ of dimension 3 and 4:
\begin{enumerate}
\item $S$ is Noetherian of global dimension 3 or 4, respectively.
\item To each $S$ one can associate an elliptic curve $E$, an automorphism $\sigma$ of $E$ and an invertible 
sheaf $\LL$ on $E$. 
\item
The following are equivalent:
\begin{itemize}
\item $S$ is a PI algebra;
\item $S$ is module finite over its center $\ZZ(S)$;
\item The automorphism $\sigma$ of $E$ has finite order.
\end{itemize}
\end{enumerate}
The centers of the 3 and 4-dimensional PI Sklyanin algebras $S$ were explicitly described in \cite{AST,ST}. 
Each such algebra $S$ was proved to be a maximal order in \cite[Theorem 2.10 and Introduction]{St}. 
Thus, we have a reduced trace map 
\begin{equation}
\label{ex5.3-tr}
\tr_{\red} : S \to \ZZ(S).
\end{equation}

It was proved in \cite{WWY1,WWY2} that each 3 and 4-dimensional PI Sklyanin algebra $S$ can be obtained as a higher order specialization, 
which gives rise to a Poisson order structure
\begin{equation}
\label{ex5.3-del}
(S, \ZZ(S), \partial_{\spp}).
\end{equation}
In \cite[Theorem 1.1]{WWY1} and \cite[Theorem 1.2]{WWY2}, it was proved that all of these Poisson order structures are nontrivial and the 
symplectic foliations of the spectrum of $\ZZ(S)$ were explicitly described.

Applying \prref{spec-red} we obtain the following:
\bpr{2}
For each 3 and 4-dimensional Sklyanin algebra $S$ corresponding to a finite order automorphism $\sigma$ of an elliptic curve $E$, 
the trace map \eqref{ex5.3-tr} and the Poisson order \eqref{ex5.3-del} are packaged into a Poisson trace order structure
\[
(S, \ZZ(S), \partial_{\spp}, \tr_{\red}).
\] 
\epr
\subsection{Cluster algebras}
\label{5.4}
Cluster algebras were defined by Fomin and Zelevinsky in \cite{FZ}. For the last 20 years they have played a 
prominent role in a number of areas of mathematics and mathematical physics. Next we show that their root of unity quantum counterparts have
Poisson trace order structures under natural assumptions. 

To define them, choose a positive integer $N$, a set of exchangeable indices $\ex\subseteq [1,N]$ and an exchange matrix $\B$, 
which is an integer matrix of size $|\ex| \times N$ 
with skew-symmetrizable principal part. Fix a positive integer $\ell$ and a 
primitive $\ell$-th root of unity $\ep^{1/2} \in \Cset$. For a skew-symmetric matrix
\[
\Om \in M_{N \times N} ( \Zset/\ell)
\]
define the mixed quantum torus/quantum affine space algebra
\[
\TT_\ep(\Om)_{\geq} := \frac{\Cset \langle x_k^{\pm 1}, x_i ; i \in \ex, k \in [1,N] \backslash \ex \rangle}{( x_i x_j - \ep^{\Om_{ij}} x_j x_i )} \cdot  
\]
(This algebra is a domain, and thus so are the rest of the algebras discussed in this subsection and Sect. \ref{6.4}, because they are constructed 
as subalgebras of $\TT_\ep(\Om)_{\geq}$.)
The matrices $\B$ and $\Om$ are called $\ell$-compatible if
\begin{equation}
\label{compat}
\ol{\B}^\top \hspace{-0.15cm} \Omega  = 
\begin{bmatrix}
\, \ol{D} \; 0
\end{bmatrix} \hspace{-0.05cm}, 
\end{equation}
where $D$ is a diagonal integer matrix with positive diagonal entries which skew-symmetrizes 
the principal part $\B$. Here, $\ol{C}$ denotes the reduction modulo $\ell$ of an integer matrix $C$. 
In this setting, one defines the root of unity quantum cluster algebra $\A_\ep(\B)$ and 
and the upper root of unity quantum cluster algebra $\U_\ep(\B)$, which satisfy 
\[
\A_\ep(\B) \subseteq \U_\ep(\B).
\]
The former algebra is the complex algebra generated by the cluster variables in all seeds obtained from the 
initial one by iterative mutations in all directions indexed by $\ex$ and the latter is the intersection of the 
mixed quantum torus/quantum affine space algebras obtained from $\TT_\ep(\Om)$ by mutation, 
see \cite[Sect. 3]{MNTY}. (In the notation for $\A_\ep(\B)$ and $\U_\ep(\B)$ one usually includes an
initial toric frame and in the construction one can allow for a subset of frozen variables to be inverted.) 

$\U_\ep(\B)$ is always a PI algebra, because it is a subalgebra of $\TT_\ep(\Om)_{\geq}$.
If it is a finitely generated algebra, then it is a maximal order, see \cite[Theorem A(1)]{HLY}, 
so we have the reduced trace map
\begin{equation}
\label{CA-red}
\tr_{\red} : \U_\ep(\B) \to \ZZ(\U_\ep(\B)).
\end{equation}

The root of unity upper quantum cluster algebra $\U_\ep(\B)$ is called {\emph{strict}} if there exists a skew-symmetric matrix
\[
\La \in M_{N \times N} (\Zset) 
\]
such that $\ol{\La} = \Om$ and the pair $(\B, \La)$ is compatible in the sense that the analog of \eqref{compat} holds over $\Zset$:
\[
\B^\top %\hspace{-0.15cm} 
\La  = 
\begin{bmatrix}
\, D \; 0
\end{bmatrix} \hspace{-0.05cm}.
\]

If 
\begin{enumerate}
\item[(i)] $\U_\ep(\B)$ is a strict  root of unity upper quantum cluster algebra,
\item[(ii)] $\ell$ is an odd positive integer that is coprime to the diagonal entries of the skew-symmetrizing matrix $D$ and
\item[(iii)] $\A_\ep(\B) =\U_\ep(\B)$, 
\end{enumerate}
then by the proof of \cite[Theorem 5.2]{MNTY}, 
$\U_\ep(\B)$ is a specialization of the corresponding upper quantum cluster algebra. This gives rise to 
the Poisson order structure
\begin{equation}
\label{U-sp}
(\U_\ep(\B), \ZZ(\U_\ep(\B)), \partial_{\spp}).
\end{equation}
By \cite[Theorem B(4)]{HLY}, if
\begin{enumerate}
\item[(iv)] $\U_\ep(\B)$ is a finitely generated algebra over $\Cset$,
\end{enumerate}
then $\U_\ep(\B)$ is a finitely generated module over $\ZZ(\U_\ep(\B))$. 

\prref{spec-red} now implies the following:
\bpr{4} For every root of unity upper quantum cluster algebra $\U_\ep(\B)$ 
satisfying assumptions (i)--(iv), the trace map \eqref{CA-red} and
the Poisson order structure \eqref{U-sp} give rise to the Poisson trace order
\[
(\U_\ep(\B), \ZZ(\U_\ep(\B)), \partial_{\spp}, \tr_{\red}).
\]
\epr

The assumptions (i)--(iv) are satisfied in wide generality. All important root of unity quantum cluster algebras are strict and 
finitely generated; (ii) is a standard assumption on $\ell$ analogous to the one in Sect. \ref{5.2}. There are a 
number of recent papers showing that cluster algebras 
equal the corresponding upper cluster algebras in many important situations 
on the classical and quantum levels \cite{CGGLSS,GY1,GY2,M}.

\subsection{Symplectic reflection algebras}\label{symprefalg}
Symplectic reflection algebras were defined and and studied as a generalisation of rational Cherednik algebras in \cite{EG}. We briefly recall the definition. Let $(V, \omega)$ be a complex symplectic vector space with basis $\{x_i : 1 \leq i \leq 2n\}$,  and let $G$ be a finite subgroup of $\mathrm{Sp}(V)$. Define an element $s$ of $G$ to be a \emph{symplectic reflection} if $\mathrm{rank}(1 - s) = 2$. The set $\mathcal{S}$ of symplectic reflections in $G$ is closed under conjugation. Choose $t \in \mathbb{C}$ and a $G$-invariant function 
$$\mathbf{c} : \mathcal{S} \longrightarrow \mathbb{C} : c \mapsto c_s.$$
For $g \in G$ there is a decomposition of $V$ as $\mathbb{C}\langle g \rangle$-module,
$$ V \; = \mathrm{im}(1 - g) \oplus \mathrm{ker}(1 - g), $$
which is easily confirmed to be $\omega$-orthogonal. For later use, note that this implies that 
\begin{equation}\label{big} \textit{if } g \neq 1 \textit{ then } \mathrm{dim}_{\mathbb{C}}(\mathrm{ker}(1 - g)) \leq 2n - 2.
\end{equation}
For $s \in \mathcal{S}$ denote by $\omega_s$  the skew-symmetric
form on $V$ which has $\mathrm{ker}(1 - s)$ as its radical and coincides with $\omega$ on $\mathrm{im}(1 - s)$. Now
define the \emph{symplectic reflection algebra} $H_{t,\mathbf{c}}$ to be the $\mathbb{C}$--algebra with generators $\{x_i,g \, : \, 1 \leq i \leq 2n, g \in G \}$ and relations those for $G$, together with
$$ gx_ig^{-1} \, = \, g(x_i), \textit{  and  } x_ix_j - x_jx_i = t\omega(x_i,x_j)1_G + \sum_{s \in \mathcal{S}} c_s\omega_s(x_i,x_j)s$$
for $i,j = 1,2, \ldots , 2n$ and $g \in G$. It is easy to see that, for $\lambda \in \mathbb{C}^{\times}$, $H_{t,\mathbf{c}} \cong H_{\lambda t ,\lambda \mathbf{c}}$, so that, as $\mathbf{c}$ varies, one need only study the algebras $H_{0,\mathbf{c}}$ and $H_{1,\mathbf{c}}$. 

It is clear from the defining relations that $H_{0,\mathbf{0}}$ is simply the skew group algebra $S(V)\ast G$. It is also obvious that an  $\mathbb{N}$-fitration $\mathcal{F}$ can be defined on $H_{t,\mathbf{c}}$ by declaring the elements of $G$ to have degree 0 and the nonzero elements of $V$ to have degree 1. It is then a fundamental theorem of Etingof and Ginzburg \cite[Theorem 1.3]{EG} that 
\begin{equation}\label{grade} \mathrm{gr}_{\mathcal{F}}H_{t,\mathbf{c}} \; \cong \; H_{0,\mathbf{0}} \; = \; S(V)\ast G. 
\end{equation}

At the level of algebras, the dichotomy between the cases $t = 0$ and $t = 1$ is stark. First, when $t = 0$, $H_{0,\mathbf{c}}$ is a finite module over its affine center which we denote by $Z_{\mathbf{c}}$, and $Z_{\mathbf{c}}$ is a deformation of the algebra of invariants $S(V)^G$, so that $\mathrm{Maxspec}(Z_{\mathbf{c}})$ is a deformation of $V/G$; these results are due to \cite[Theorem 3.1]{EG}. But when $t = 1$, $Z(H_{1,\mathbf{c}}) \cong \mathbb{C}$ by \cite[Proposition 7.2(2)]{BG}. For all values of the parameters $(t,\mathbf{c})$, it follows routinely from filtered-graded technology that $H_{t,\mathbf{c}}$ is a prime Noetherian algebra with good homological properties.

By means of the specialisation mechanism of \prref{spec-red}(1), $H_{0,\mathbf{c}}$ is a Poisson $Z_{\mathbf{c}}$-order; 
the details are explained in \cite[$\S$2.2]{BG}. Moreover, it is shown in \cite[Theorem 7.8]{BG} that 
$\maxspec(Z_{\mathbf{c}})$ has only finitely many symplectic leaves, and that these coincide 
with the irreducible components of the rank stratification of the Poisson bracket. 

To see that $H_{0,\mathbf{c}}$ is a Poisson trace order with respect to the reduced trace, we need the following result from \cite[Theorem 4.4]{B0}. For the convenience of the reader we give the proof here.

\bpr{symp-max} Retain the notation and hypotheses introduced above. Then $H_{0, \mathbf{c}}$ is a maximal order, and hence its center 
$Z_{\mathbf{c}}$ is an integrally closed domain.
\epr 

\begin{proof} Thanks to (\ref{big}), if $1 \neq g \in G$ then the ideal $\langle v - v^g : v \in V \rangle$ of the polynomial algebra $S(V)$ is not contained in any height one prime of $S(V)$. The skew group algebra $H_{0,\mathbf{0}} = S(V)\ast G$ is therefore a maximal order by \cite[Theorem 4.6]{M}. In view of (\ref{grade}) it follows from \cite[Theorem 5]{V2} that we can lift the maximal order property to $H_{0,\mathbf{c}}$ as required. The second claim follows since the centre of a prime Noetherian maximal order is an integrally closed domain by  \cite[Proposition 5.1.10(b)(i)]{McR}.
\end{proof}

The following result is now immediate from Propositions \ref{pspec-red} and \ref{psymp-max}.

\bpr{symp-red} For any choice of data $(V,\omega)$, $G$ and $\mathbf{c}$ as explained above, 
$(H_{0,\mathbf{c}}, Z_{\mathbf{c}}, \partial_{\spp}, \tr_{\red})$ is a Poisson trace order.
\epr
%%%%%%%%%%%
\sectionnew{Construction of Poisson trace orders: partial centers}
\label{6}
In this section we use Theorems \ref{treg-trace1} and \ref{tbase-change} 
to construct large classes of Poisson trace orders with respect to proper central subalgebras
from the areas of quantum groups, Nichols algebras and cluster algebras. 

\subsection{Quantized enveloping algebras and quantum Schubert cell algebras at roots of unity}
\label{6.1}
Consider the families (1) and (3) in Sect. \ref{5.2}. De Concini, Kac and Procesi constructed central subalgebras of 
the quantized universal enveloping algebras $\UU_\ep(g)$ and the quantum Schubert cell algebras $\UU_\ep^w$. The central subalgebra 
\[
Z_0(\UU_\ep(\g))
\]
of $\UU_\ep(g)$ was constructed in \cite[Corollary 3.1]{DK} as the subalgebra generated by the $\ell$-th powers of the root vectors $E_\alpha$, $F_\alpha$ 
and Chevalley generators $K_i^{\pm1}$ of $\UU_\ep(g)$. Recall that $\ell$ denotes the order of the root of unity $\ep \in \Cset$. The central subalgebra 
\[
Z_0(\UU_\ep^w)
\]
of $\UU_\ep^w$ was constructed in \cite[Proposition 3.1]{DKP2} as the subalgebra generated by the $\ell$-th powers of the root vectors $E_\alpha$
of $\UU_\ep^w$. 

In \cite{DKP1,DKP2} it was proved that $Z_0(\UU_\ep(g))$ and $Z_0(\UU_\ep^w)$ are Poisson subalgebras 
of the full centers $\ZZ(\UU_\ep(g))$ and $\ZZ(\UU_\ep^w)$ with respect to the underlying Poisson algebra structures 
of the first and third Poisson orders in \eqref{PO}. Therefore, those two Poisson orders restrict to the Poisson orders
\begin{equation}
\label{ex6.1-del}
(\UU_\ep(\g), Z_0(\UU_\ep(\g)), \partial_{\spp}|_{Z_0(\UU_\ep(\g))}) \quad \mbox{and} \quad 
(\UU_\ep^w, Z_0(\UU_\ep^w), \partial_{\spp}|_{Z_0(\UU_\ep^w)}).
\end{equation}
Furthermore, because of PBW bases, $\UU_\ep(\g)$ is a free module over $Z_0(\UU_\ep(\g))$ of rank $\ell^{\dim \g}$ and $\UU_\ep^w$ 
is a free module over $Z_0(\UU_\ep^w)$ of rank $\ell^{l(w)}$, where $l(w)$ is the length of $w$. 
Thus, we have the regular trace maps
\begin{equation}
\label{ex6.1-tr}
\tr_{\reg} : \UU_\ep(\g) \to Z_0(\UU_\ep(\g)) \quad \mbox{and} \quad
\tr_{\reg} : \UU_\ep^w \to Z_0(\UU_\ep^w).
\end{equation}
\thref{reg-trace1} in combination with \eqref{ex6.1-del} and \eqref{ex6.1-tr} implies the following:
\bpr{b1} For all complex simple Lie algebras $\g$, 
Weyl group elements $w$ and primitive roots of unity $\ep \in \Cset$ of odd order $\ell$ such that $3 \nmid \ell$ if $\g$ is of type $G_2$,
the quadruples
\[
(\UU_\ep(\g), Z_0(\UU_\ep(\g)), \partial_{\spp}|_{Z_0(\UU_\ep(\g))}, \tr_{\reg}) \quad \mbox{and} \quad 
(\UU_\ep^w, Z_0(\UU_\ep^w), \partial_{\spp}|_{Z_0(\UU_\ep^w)}, \tr_{\reg}).
\]
are Poisson trace orders.
\epr 
This result also follows by applying \thref{reg-trace1} and then doing a base change via \thref{base-change} because 
the regular trace map $\tr_{\reg} : \UU_\ep(\g) \to Z_0(\UU_\ep(\g))$ equals the composition of 
\[
\tr_{\red} : \UU_\ep(\g) \to \ZZ(\UU_\ep(\g)) \quad \mbox{and} 
\quad \tr_{\ZZ(\UU_\ep(\g))/Z_0(\UU_\ep(\g))} : \ZZ(\UU_\ep(\g)) \to Z_0(\UU_\ep(\g))
\]
and the regular trace map $\tr_{\reg} : \UU_\ep^w \to Z_0(\UU_\ep^w)$ equals the composition of
\[
\tr_{\red} : \UU_\ep^w \to \ZZ(\UU_\ep^w) \quad \mbox{and} \quad 
\tr_{\ZZ(\UU_\ep^w)/Z_0(\UU_\ep^w)} : \ZZ(\UU_\ep^w) \to Z_0(\UU_\ep^w).
\]
\thref{base-change} is applicable because $Z_0(\UU_\ep(\g))$ and $Z_0(\UU_\ep^w)$ are normal
(since they are tensor products of polynomial and Laurent polynomial algebras) and 
$\ZZ(\UU_\ep(\g))$ and $\ZZ(\UU_\ep^w)$ are finitely generated modules over them as 
submodules of the finitely generated modules $\UU_\ep(\g)$ and $\UU_\ep^w$.

\subsection{Quantized function algebras at roots of unity}
\label{6.2}
Consider the family (2) in Sect. \ref{5.2} of the quantized coordinate ring at root of unity $F_\ep[G]$ of 
a connected, simply connected complex simple algebraic group $G$.  
De Concini and Lyubashenko constructed in \cite[Proposition 6.4]{DL} a 
central subalgebra $F_0[G]$ of $F_\ep[G]$ such that
\[
F_0[G] \cong \Qset(\ep)[G].
\]
It was proved in \cite[Proposition 2.2]{BG0} and \cite[Theorem on p. 1]{BGS} that $F_\ep[G]$ is a 
free module over $F_0[G]$ of rank $\ell^{\dim G}$. This gives rise to the regular trace maps 
\begin{equation}
\label{tr-OG}
\tr_{\reg} : F_\ep[G] \to F_0[G].
\end{equation}

In addition, in \cite[Sect. 8]{DL} it was proved that $F_0[G]$ is a Poisson subalgebra of $\ZZ(F_\ep[G])$ with respect to the 
underlying Poisson algebra structure of the second Poisson order in \eqref{PO}, and that it is isomorphic to the standard Poisson 
algebra structure on the coordinate ring of $G$ over $\Qset(\ep)$. Hence, by restriction we get the Poisson order
\begin{equation}
\label{del-OG}
(F_\ep[G], F_0[G], \partial_{\spp}|_{F_0[G]}).
\end{equation}
Applying \thref{reg-trace1}, we obtain:
\bpr{b3} 
For all connected simply connected complex simple algebraic groups $G$, and primitive roots of unity $\ep \in \Cset$ of odd order $\ell$ such that $3 \nmid \ell$ if $\g$ is of type $G_2$,
the trace map \eqref{tr-OG} and Poisson order structure \eqref{del-OG} give rise to the Poisson trace order
\[
(F_\ep[G], F_0[G], \partial_{\spp}|_{F_0[G]}, \tr_{\reg}).
\]
\epr

This result can be also derived by first applying \thref{reg-trace1} and then doing a base change via \thref{base-change} because 
the regular trace map $\tr_{\reg} : F_\ep[G] \to F_0[G]$ equals the composition of 
\[
\tr_{\red} : F_\ep[G] \to \ZZ(F_0[G])  \quad \mbox{and} 
\quad \tr_{\ZZ(F_\ep[G])/F_0[G]} : \ZZ(F_\ep[G]) \to F_0[G]. 
\]
\thref{base-change} is applicable because $F_0[G]$ is normal
(since it is isomorphic to the coordinate ring of $G$ over $\Qset(\ep)$) and 
$\ZZ(F_\ep[G])$ is a finitely generated module over it, being a submodule 
of the finitely generated module $F_\ep[G]$.
\subsection{Nichols algebras of diagonal type}
\label{6.3}
In \cite{AAY} Poisson order structures were constructed on the quantum doubles of the
bozonizations of all distinguished pre-Nichols algebras of diagonal type \cite{A} with finite root systems that belong 
to a one-parameter family. These algebras play an important role in the Andruskiewitsch--Schneider program \cite{AS} on 
classifying finite dimensional pointed Hopf algebras; in the diagonal case such are classified via
Heckenberger's Weyl groupoid \cite{He}. The algebras in this family are over $\Cset$ and are denoted by 
\[
U_{\mathfrak{q}}, \; \;  \mbox{where ${\mathfrak{q}} = (q_{ij}) \in M_{r \times r} (\Cset)$ is a braiding matrix}. 
\] 
All contragredient quantum super groups at roots of unity arise as special members in this family. 
Each algebra $U_{\mathfrak{q}}$ contains a central subalgebra $Z_{\mathfrak{q}}$ 
constructed in \cite[Sect. 4.5]{AAY}. In \cite[Theorem A(a)]{AAY}, a Poisson order structure
\begin{equation}
\label{del-Nich}
(U_{\mathfrak{q}}, Z_{\mathfrak{q}}, \partial_{\spp})
\end{equation}
was constructed using specialization and Hopf theoretic methods. This structure was shown to be nontrivial and the underlying 
Poisson algebra $Z_{\mathfrak{q}}$ was proved to be isomorphic to the coordinate ring of the dual Poisson Lie group 
to an associated complex reductive Lie group \cite[Theorem A(b)]{AAY}. 
The algebra $U_{\mathfrak{q}}$ 
is a free $Z_{\mathfrak{q}}$-module of finite rank \cite[Sect. 1.3.3]{AAY}, so we have the regular trace map
\begin{equation}
\label{tr-Nich}
\tr_{\reg} : U_{\mathfrak{q}} \to Z_{\mathfrak{q}}. 
\end{equation}
\thref{reg-trace1} now implies the following:
\bpr{b3} For the quantum doubles $U_{\mathfrak{q}}$ of the bozonizations of all distinguished pre-Nichols algebras of diagonal type 
with finite root systems that belong to a one-parameter family, the Poisson order structure \eqref{del-Nich} and the trace map \eqref{tr-Nich} 
give rise to the Poisson trace orders 
\[
(U_{\mathfrak{q}}, Z_{\mathfrak{q}}, \partial_{\spp}, \tr_{\reg}).
\]
\epr
Poisson trace orders for the underlying Nichols algebras $U_{\mathfrak{q}}^+$ and their bosonizations 
$U_{\mathfrak{q}}^{\geq}$ can be similarly constructed using \thref{reg-trace1} and \cite[Theorems B and D]{AAY}. 
\subsection{Cluster algebras}
\label{6.4}
Assume the setting of Sect. \ref{5.4}. We have the central subalgebra
\[
\LL_\ep( \Om)_{\geq}:= \Cset [ (x_1)^{\pm \ell}, \ x_i^\ell  ; 
k \in \ex, i \in [1,N] \backslash \ex ] \subset \ZZ( \TT_\ep(\Om)_\geq)
\]
giving rise to the central subalgebra $\C\U_\ep(\B)$ of $\U_\ep(\B)$ obtained by 
intersecting the analogs of $\LL_\ep( \Om)_{\geq}$ over all seeds obtained from 
the initial one by mutation, see \cite[Sect. 3.1]{HLY}. By \cite[Proposition 3.1]{HLY}, if 
\begin{equation}
\label{coprime}
\mbox{$\ell$ {\it is odd and coprime to the diagonal entries of} $D$}, 
\end{equation}
cf. \eqref{compat}, then 
\[
\C\U_\ep(\B) \cong \U(\B),
\] 
where $\U(\B)$ is the upper cluster algebra associated to the exchange matrix $\B$.

The algebra $\TT_\ep( \Om)_{\geq}$ is a free module over $\LL_\ep( \Om)_{\geq}$ of rank $\ell^N$, so 
we have the regular trace map
\[
\tr_{\reg} : \TT_\ep( \Om)_{\geq} \to \LL_\ep( \Om)_{\geq}. 
\]
By \cite[Theorem B(1)--(3)]{HLY}, if \eqref{coprime} holds, then $\tr_{\reg}$ restricts to a trace map
\begin{equation}
\label{tr-U-CU}
\tr_{\reg}|_{\U_\ep(\B)} : \U_\ep(\B) \to \C\U_\ep(\B)
\end{equation}
making the triple $(\U_\ep(\B), \C\U_\ep(\B), \tr_{\reg}|_{\U_\ep(\B)})$ a Cayley--Hamilton algebra of degree $\ell^N$. 
We note that $\U_\ep(\B)$ is very rarely a free $\C\U_\ep(\B)$-module, so, here \eqref{tr-U-CU} is 
not a regular trace for the pair $(\U_\ep(\B), \C\U_\ep(\B))$ and \thref{reg-trace1} is not applicable. 

If conditions (i)--(iv) in Sect. \ref{5.4} are satisfied then the trace map \eqref{tr-U-CU} equals the composition of 
\[
\tr_{\red} :  \U_\ep(\B) \to \ZZ( \U_\ep(\B) ) 
\quad \mbox{and} 
\quad \tr_{ \ZZ(\U_\ep(\B)) / \C\U_\ep(\B)} :  \ZZ(\U_\ep(\B)) \to  \C\U_\ep(\B). 
\]
The last trace is defined because $\C\U_\ep(\B) \cong \U(\B)$ is a normal algebra since the latter algebra is 
an intersection of a collection of algebras which are tensor products of polynomial and Laurent polynomial algebras, 
see \cite[Sect. 3.1]{HLY}. 

By \cite[Theorem 5.2]{MNTY}, under the assumptions (i)--(iii), $\C\U_\ep(\B)$ is a Poisson subalgebra of $\ZZ(\U_\ep(\B))$ 
under the underlying Poisson algebra structure of the Poisson order in \eqref{U-sp}. This gives rise to the Poisson order structure
\[
(\U_\ep(\B), \C\U_\ep(\B), \partial_{\spp}|_{\C\U_\ep(\B)}).
\] 
\thref{base-change} implies that 
\bpr{b4} For every root of unity quantum cluster algebra $\U_\ep(\B)$ satisfying assumptions 
(i)--(iv) in Sect. \ref{5.4}, the quadruple
\[
(\U_\ep(\B), \C\U_\ep(\B), \partial_{\spp}|_{\C\U_\ep(\B)}, \tr_{\reg}|_{\U_\ep(\B)})
\]
is a Poisson trace order.
\epr
%%%%%%%%%%%%%%%%%%%%%% References %%%%%%%%%%%%%%%%%%%%%%%%%%%%%%%%%%%%%%%

%%%%%%%%%%%%%%%%%%%%%%%%%%%%%%%%%%%%%%%%%%%%%%%%%%%%%%%%%%%%%%%%%%%%%%%%%%%%%%%
%%%%%%%%%%%%%%%%%%%%%%%%%%%%%%%%%%%%%%%%%%%%%%%%%%%%%%%%%%%%%%%%%%%%%%%%%%%%%%
\end{document}